\documentclass[11pt,equ]{amsart}

\usepackage{graphicx,pdfsync}
\usepackage{hyperref}

\newcommand{\eqn}[1]{(\ref{#1})}
\newcommand{\sA}{{\mathcal A}}

\newcommand{\df}{\displaystyle\frac}

\newcommand{\hsp}{\hspace{\parindent}}

\newcommand{\In}{\infty}

\newcommand{\CC}{{\Bbb C}}

\newcommand{\RR}{{\Bbb R}}

\newcommand{\ZZ}{{\Bbb Z}}

\newcommand{\beql}[1]{\begin{equation}\label{#1}}
\newcommand{\eeq}{\end{equation}}

\begin{document}

{\bf Reference.} {\em  The Ultimate Challenge: The $3x+1$ Problem}. Edited by Jeffrey C. Lagarias. American Mathematical Society, Providence, RI, 2010, pp. 3--29. 
\bigskip
\bigskip
\bigskip
\title{The $3x+1$ Problem: An   Overview}
\author{Jeffrey C. Lagarias}
\address{Department of Mathematics,
University of Michigan,
Ann Arbor, MI 48109-1043}
\email{lagarias@umich.edu} 


\maketitle

%
%
%

\section{Introduction}

The $3x+1$ problem concerns 
the following  innocent seeming arithmetic procedure applied to  integers: If an integer $x$  is odd
then ``multiply  by three and add one", while if it is even then  ``divide by two".
This operation is described by the  {\em Collatz function}
$$
C(x) = 
\left\{
\begin{array}{cl}
3x+1 & \mbox{if}~ x \equiv 1~~ (\bmod ~2 ) , \\
~~~ \\
\df{x}{2} & \mbox{if} ~~x \equiv 0~~ (\bmod~2) .
\end{array}
\right.
$$
The $3x+1$ problem, which is often called the 
{\em Collatz problem},  concerns  the behavior of this function under iteration,
starting with a given positive integer $n$.  

\smallskip
{\bf  $3x+1$ Conjecture.}
{\em Starting from any positive integer $n$,  
iterations of the function $C(x)$  will eventually  reach the number $1$.
Thereafter iterations will cycle, taking successive values $1, 4, 2, 1,...$.}

\smallskip
This problem  goes under many other names, including
the  {\em Syracuse problem}, {\em Hasse's algorithm},
 {\em Kakutani's problem} and {\em Ulam's problem}.  

\smallskip
 A commonly used reformulation of the $3x+1$  problem iterates a different function,
the {\em $3x+1$ function}, given by
$$
T(x) = 
\left\{
\begin{array}{cl}
\df{3x+1}{2} & \mbox{if}~ x \equiv 1~~ (\bmod ~2 ) , \\
~~~ \\
\df{x}{2} & \mbox{if} ~~x \equiv 0~~ (\bmod~2) .
\end{array}
\right.
$$
From the viewpoint of iteration the two functions are simply related;
iteration of $T(x)$ simply omits some steps in the iteration of the
Collatz function $C(x)$. 
The relation of  the $3x+1$ function $T(x)$ to  the Collatz function $C(x)$  is that:
$$
T(x) = 
\left\{
\begin{array}{cl}
C( C(x)) & \mbox{if}~ x \equiv 1~~ (\bmod ~2 ) ~, \\
~~~ \\
C(x) & \mbox{if} ~~x \equiv 0~~ (\bmod~2) ~.
\end{array}
\right.
$$
As it turns out, the function $T(x)$ proves more convenient
 for analysis of the problem
 in a number of significant ways, as first
  observed independently by  Riho Terras (\cite{Ter76}, \cite{Ter79}) and by C. J. Everett \cite{Ev77}.

\smallskip
The $3x+1$  problem  has fascinated
mathematicians and non-mathematicians alike. 
It has been studied by mathematicians, physicists, and computer scientists. 
It remains an unsolved problem, which appears to be extremely difficult.

\smallskip
This paper aims to address two questions: \smallskip
 %
 %
\begin{enumerate}
\item[(1)]
{\em What can mathematics  currently say about this problem?}\medskip

\item[(2)]
{\em  How can  this problem be hard, when it is so easy to state? }
\end{enumerate}
 %
 %

\smallskip
To address the first question, 
this overview discusses the history of work on
the problem. Then it describes generalizations of
the problem, and lists the different
fields of mathematics on which the
problem impinges.  It gives a brief summary of the current
strongest results on the problem. 

\smallskip
Besides the results summarized  here,
this  volume contains more detailed surveys  
of mathematicians' understanding of the
$3x+1$ problem and  its generalizations. 
These cover   both 
rigorously proved  results and heuristic predictions made
using probabilistic models.  
The book includes  several survey articles, 
it reprints several early papers on the problem, with commentary, 
and it   presents  
an annotated bibliography of work
on the problem and its generalizations. 

\smallskip
 To address the second question,  let us remark first
 that the true level of difficulty of any problem can only be 
determined when (and if)  it is solved. Thus there can be 
no definitive answer regarding  its difficulty.
The track record on the $3x+1$ problem so far 
 suggests that this is an extraordinarily
difficult problem, completely out of reach of present day
mathematics.  
Here we will only say that part of the  difficulty appears
to reside in an inability to analyze the pseudorandom nature of successive iterates
of $T(x)$,
which could conceivably encode very difficult computational problems.
We elaborate on this answer in \S7.

\smallskip
Is the $3x+1$  problem an important problem? 
Perhaps not  for its individual sake, where it merely stands
as a challenge. It seems to be a prototypical example of an
extremely simple to state, extremely hard to solve, problem.
A middle of the road viewpoint is that this problem is representative 
of a large class of problems,  concerning the behavior under iteration of maps that are expanding on
part of their domain and contracting on another part of their domain.
This general class of problems is  of definite importance, and is 
currently of great interest as an area
of mathematical (and physical) research; for some
perspective, see Hasselblatt and Katok \cite{HK95}. 
Progress on general methods of
solution for functions in this class would be extremely significant.

\smallskip
This overview describes where things currently stand on the $3x+1$ problem
and how it relates to various fields of mathematics. For a detailed
introduction to the problem, 
see the following paper of Lagarias \cite{Lag85} (in this volume).
In \S2 we give some history of the problem; this presents some new information beyond
that given in  \cite{Lag85} .
Then in \S3 we give a flavor of the behavior of the $3x+1$ iteration.
 In \S4 we discuss various frameworks for generalizing
the problem; typically these concern iterations of functions having a similar
appearance to the $3x+1$ function. In \S5 we review areas of research:  
these comprise different
fields of mathematics and computer science on which this problem impinges.
In \S6 we summarize the current best results on the problem in various directions.
In \S7 we discuss the hardness of the $3x+1$ problem.
In \S8 we describe some research directions for future progress.
In \S9 we address the question: ``Is the $3x+1$ problem a good problem?"
In the concluding section \S10 we offer some advice on working on
$3x+1$-related problems.

%
%
%

\section{History and Background}

The  $3x+1$ problem circulated by word
of mouth for many years.
It is generally attributed to Lothar Collatz. He has stated (\cite{Col80}) that he
took lecture courses  in 1929 with Edmund Landau   and Fritz von Lettenmeyer in
G\"{o}ttingen, and courses in 1930 with Oskar Perron in Munich and with Issai Schur
in Berlin, the latter course including some graph theory.
  He was interested in graphical representations of iteration of  functions. 
In his notebooks in the 1930's he formulated questions 
on iteration of arithmetic functions of a similar kind  (cf. \cite[p.\,3]{Lag85}). 
Collatz  is said by others to have circulated the problem orally at the
International Congress of Mathematicians in Cambridge, Mass.
 in 1950. Several  people whose names were subsequently associated with the problem
gave invited talks at this International Congress, including H. S. M. Coxeter, S. Kakutani, 
and S. Ulam. 
 Collatz \cite{Col86} (in this volume)   states that  he described the $3x+1$ problem
to Helmut Hasse in 1952 when they were colleagues at the University of Hamburg.
Hasse was  interested in the problem, and  wrote 
about it in lecture notes  in 1975 (\cite{Ha75}). 
Another claimant  to having originated the $3x+1$ problem is  Bryan  Thwaites \cite{Th85},
who asserts that he came up with the problem in 1952. Whatever is its true
origin, the
 $3x+1$ problem was already
 circulating at the University of Cambridge in the late
 1950's, according to John H. Conway and to Richard Guy \cite{Guy09}.

\smallskip
 There was no published mathematical literature about the $3x+1$ problem until
 the early 1970's. This may have been, in part, because the 1960's was 
a period dominated by  Bourbaki-style mathematics. The
Bourbaki viewpoint  emphasized complete presentations of 
theories with rich internal structure,  which interconnect
with other areas of core mathematics (see Mashaal \cite{Mas06}).   In contrast, 
the $3x+1$ problem initially appears to be an isolated problem unrelated
to the rest of mathematics. 
Another obstacle was the difficulty in proving interesting results
about the $3x+1$ iteration.
The results that could be proved 
appeared pathetically weak, so that it could seem  damaging to one's professional
reputation to publish them. In some mathematical circles it might
have seemed in bad taste  even to show interest in
such a problem, which appears d\'{e}class\'{e}.

\smallskip
During the 1960's, various problems related to the $3x+1$ problem
appeared in print, typically as unsolved problems. This included one
of the original problems of Collatz from the 1930's, which concerned the
behavior under iteration of the function 
$$U(2n)= 3n, ~U(4n+1)= 3n+1, ~U(4n+3)= 3n+2.
$$
The function $U(n)$ defines a permutation of the integers, and the question concerns whether
the iterates of the value $n=8$ form an infinite set.  This problem was
raised by Murray Klamkin  \cite{Klm63} in 1963 (see  Lagarias \cite[p.\, 3]{Lag85}),
and remains unsolved.
 Another such problem was 
posed by Ramond Queneau, a founder of the French mathematical-literary
group Oulipo (Ouvroir de litt\'{e}rature potentielle), which concerns
allowable rhyming patterns generalizing those used in poems by the 12-th
century troubadour, Arnaut Daniel. This problem turns out to be  related
 to a $(3x+1)$-like
function whose behavior under iteration is exactly analyzable, see Roubaud \cite{Rou69}.
Concerning the $3x+1$ problem itself, during the 1960's  
 large computations were done  testing the truth
of the conjecture. These reportedly verified the conjecture for all $n \le 10^{9}$. 

\smallskip
To my knowledge,  the $3x+1$ problem first appeared  in print in 1971, in the written version
of a 1970 lecture  by H. S. M. Coxeter \cite{Cox71} (in this volume). It was presented there
``as a piece of mathematical gossip."
In 1972 it appeared in six different publications, including a 
Scientific American column by Martin Gardner \cite{Gar72}
that gave it wide
publicity. Since then
there has been  a steady stream of work on it, 
now amounting to several hundred  publications. 

\smallskip
Stanislaw Ulam  was one of many who circulated the problem;
the name ``Ulam's problem" has been attached to it in some 
circles.  He was a pioneer in ergodic theory and very interested in  iteration of functions
and their study by computer; he formulated many problem lists (e.g. \cite{Ulam64}, \cite{Ulam-CC}). 
A collaborator, Paul Stein \cite[p. 104]{Ste89}, wrote about Ulam:

\begin{quotation}
Stan was not a number theorist, but he knew many number-theoretical facts.
As all who knew him well will remember, it was Stan's particular pleasure
to pose difficult, though simply stated, questions in many branches of
mathematics. Number theory is a field particularly vulnerable to the
``Ulam treatment," and Stan proposed more than his share of hard questions;
not being a professional in the field, he was under no obligation to
answer them.
\end{quotation}

\noindent Ulam's  long term collaborator 
C. J. Everett \cite{Ev77} wrote one of the early papers about the 
$3x+1$ problem in 1977.

\smallskip
The $3x+1$ problem can also be formulated in the backwards direction,
as that of determining the smallest set $S_0$ of integers containing $1$   
which is closed under the affine maps $x \mapsto 2x$ and $3x+2 \mapsto 2x+1$,
where the latter map may only be applied to inputs $3x+2$ whose output $2x+1$ will be an
integer. The $3x+1$ conjecture then asserts that $S_0$ will be the set of all positive integers.
This connects the $3x+1$ problem with problems  on sets of integers which are 
closed under the action of affine maps.  Problems of this sort were raised by 
Isard and Zwicky \cite{IZ70} in 1970. In 1970-1971  David Klarner began
studying sets of integers closed under iteration of affine maps, 
leading to joint work with Richard Rado \cite{KR74},
published in 1974.   Interaction of
Klarner and Paul Erd\H{o}s at the University of Reading in 1971  led to  
the formulation of a (solved) Erd\H{o}s prize problem:
 Does the smallest set $S_1$ of integers containing $1$
and closed under the affine maps $x \mapsto 2x+1, x \mapsto 3x+1$ and $x \mapsto 6x+1$
have a positive (lower asymptotic) density?
This set $S_1$ was  
proved to have zero density   by D. J. Crampin and A. J. W. Hilton (unpublished),
according to  Klarner \cite{Kla82}. 
The solvers  collected $\pounds 10$ from Erd\H{o}s (\cite{Hil10}). Later Klarner \cite[p. 47]{Kla82}
formulated a revised problem: 

\smallskip
{\bf Klarner's Integer Sequence Problem. } {\em Does the smallest set of integers $S_2$ containing $1$
and closed under the affine maps $x \mapsto 2x, x \mapsto 3x+2$ and $x \mapsto 6x+3$
have a positive (lower asymptotic) density?}  

\smallskip
\noindent This problem remains
unsolved; see the paper of Guy \cite{Guy83} (in this volume) and 
accompanying editorial commentary.

\smallskip
Much  early work on the problem appeared in unusual places, some of it
in technical reports, some in problem journals. 
The annotated bibliography given in this book \cite{Lag-B1} covers some of this 
literature, see also its sequel \cite{Lag-B2}.  Although the problem began life as a curiosity, 
its general connection with various other areas of 
mathematics, including number theory, dynamical systems and theory of computation,
have made it a respectable topic for mathematical research.
A number of  very well known mathematicians have contributed results on it, including
John H. Conway \cite{Con72} and Yakov G. Sinai \cite{Sin03a}, \cite{Sin03b}.

%
%
%
\section{$3x+1$ Sampler}

 The fascination of the $3x+1$ problem involves its simple definition
 and the apparent complexity of its behavior under iteration: there
 seems to be  no simple relation between the input value $n$ and the
 iterates of $n$.   
 Exploration of its structure
 has led to the formulation of a web of subsidiary conjectures about the behavior
of iterates of the $3x+1$ function and  generalizations; these include 
conjectures (C1)--(C5) listed in \S8. Many of these conjectures seem 
to be extremely difficult problems as well,  and their exploration
has led to much further research.
Since other papers in this volume give much more information on
 this complexity, here we give only a brief sampler of  $3x+1$ function  behavior.

 %
 %
 
 \subsection{Plots of Trajectories } 
 
\begin{figure}[h]
\centering
\includegraphics[width=3.5in]{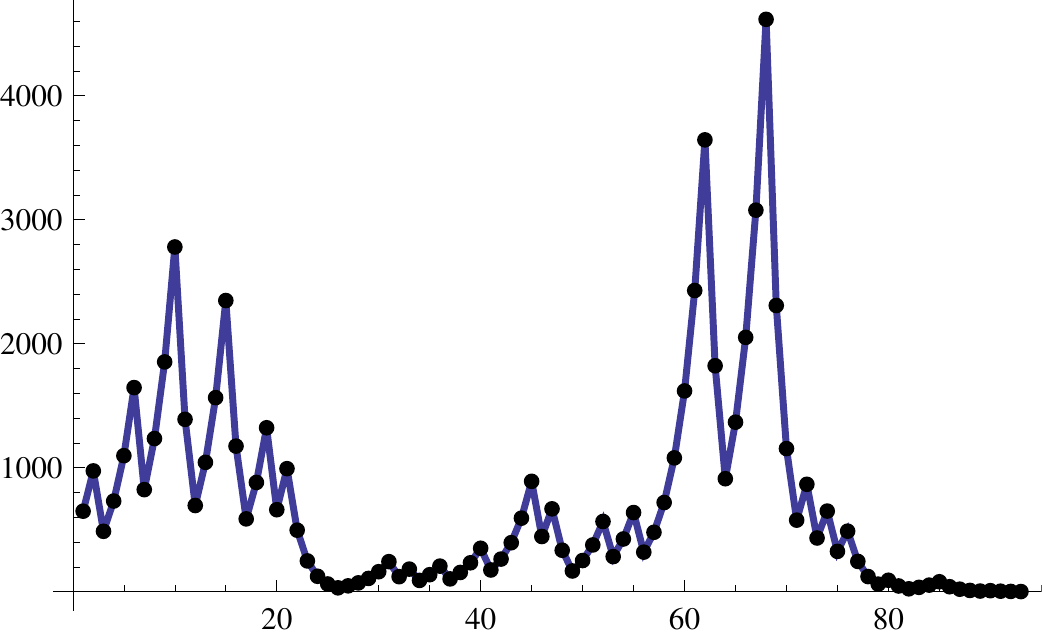}
\caption{Trajectory of $n=649$ plotted on standard vertical scale}
\label{fig21}
\end{figure}

 By the {\em trajectory} of $x$ under a function $T$, we mean the  forward orbit of $x$,
 that is, the  sequence of its forward iterates \\
 $( x, T(x), T^{(2)}(x), T^{(3)}(x), ...)$. 
 Figure \ref{fig21}  displays the $3x+1$-function iterates of $n=649$ plotted on a standard scale.
 We see an irregular series of  increases and decreases, leading to the
  name ``hailstone numbers"  proposed by Hayes \cite{Hay84},
 as hailstones form by repeated upward and downward movements
 in a thunderhead.
 
%
%
%

\begin{figure}
\centering
\includegraphics[width=3.5in]{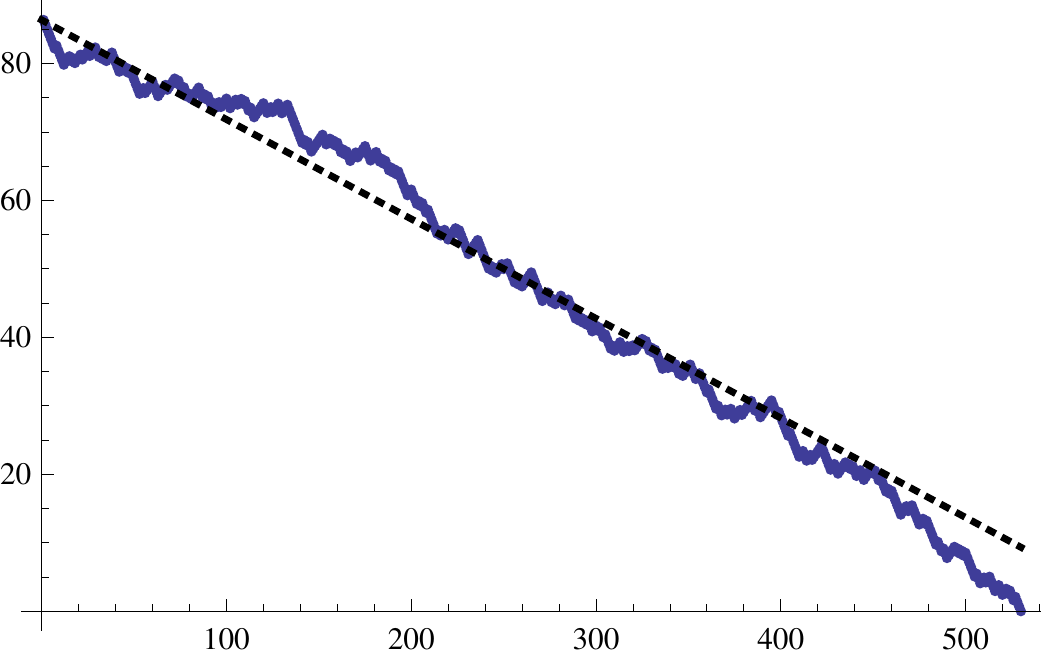}
\caption{Trajectory of $n_0=100 \lfloor \pi \cdot 10^{35}\rfloor$
 plotted on  a logarithmic vertical  scale $y= \log T^{(x)}(n_0)$ (natural logarithm), $x$ an  integer. The dotted line is a
 probability model prediction for a ``random" trajectory for this size $N$.}
\label{fig22}
\end{figure}

\smallskip
To gain insight into  a problem it helps to choose an appropriate scale for
picturing it. Here it  is useful to view long trajectories on a
 logarithmic scale, i.e., to plot $\log T^{(k)}(n)$ versus $k$. 
 Figure \ref{fig22} displays the iterates of 
 $ n_0 = 100\lfloor \pi 10^{35} \rfloor $
  on such a scale. Using this scale we see 
  a decrease at a certain geometric rate to the value of $1$, indicated by the
 trajectory having roughly a constant slope. 
 This is characteristic of most long trajectories. As explained in \S3.3
 a probabilistic model predicts that
 most trajectories plotted on a logarithmic scale will stay  close  to a line of constant
 slope  $-\frac{1}{2} \log \frac{3}{4} \sim -0.14384,$
 thus taking about $6.95212 \log n$ steps to reach $1$. This line is
 pictured as the dotted line in Figure \ref{fig22}. This trajectory
 takes $529$ steps to reach $n=1$, while the probabilistic model
 predicts about $600$ steps will be taken.

%
%
%

 \smallskip
 On the other hand, plots of trajectories suggest that
  iterations of the $3x+1$ function also seem to
exhibit pseudo-random features, i.e.  the
successive iterates of a random starting value
seem to increase or decrease in an unpredictable manner.
 From this perspective there are some regularities of the iteration 
that appear (only) describable  as  statistical in nature: they are
assertions about the majority of trajectories in ensembles of trajectories rather than about
individual trajectories. \\

 %
 %

\subsection{Patterns} 
Close examination of  the iterates
of  the 
$3x+1$  function $T(x)$ for different starting values reveals a myriad of internal patterns. 
A simple pattern is that the initial iterates of 
$n=2^m-1$ are
$$
T^{(k)}(2^m-1) = 3^{k} \cdot 2^{m-k} -1, ~~\mbox{for} ~1 \le k \le m.
$$
In particular,  $T^{(m)}(2^m -1) = 3^m -1$; this 
example  shows that the iteration can sometimes
reach  values  arbitrarily larger than the initial
value, either on an absolute or a relative scale, 
even if, as conjectured, the iterates eventually reach $1$.
Other patterns include the appearance of
occasional  large clusters of consecutive numbers  which all take exactly the
same number of iterations to reach the value $1$. Some of these patterns are
easy to analyze, others are more elusive.

\smallskip
Table 2.1 presents   data on iterates of the $3x+1$ function $T(x)$ for $n= n_0+ m$, 
$0 \le m = 10j +k \le 99$, with 
$$
n_0 = 100\lfloor \pi \cdot 10^{35} \rfloor =
31,415,926,535,897,932,384,626,433,832,795,028,800.
$$
Here $\sigma_{\infty} (n)$ denotes the {\em total stopping time} for $n$, which counts
the number of iterates of the $3x+1$-function $T(x)$ needed to reach $1$ starting from $n$,
counting $n$ as the $0$-th iterate. This number is the same as  the number
of even numbers appearing in the trajectory of the Collatz function before first reaching $1$.

\medskip
\begin{center}
\begin{tabular}{|c|c|c|c|c|c|c|c|c|c|c|}
 \hline
 & $j=0$ & $j=1$ & $j=2$ & $j=3$ & $j=4$ & $j=5$ & $j=6$ & $j=7$ & $j=8$ & $j=9$ \\ \hline
$k=0$& 529&  529& 529 & 678& 529 & 529& 846 & 529& 846& 846\\
$k=1$ & 659&  659& 529 & 678& 659 & 529 & 846& 529& 529& 529\\
$k=2$ & 846&  529& 659 & 529& 529 & 529& 659& 846& 529& 659\\
$k=3$ & 846&  529& 659 & 846 & 659 & 529 & 659& 846 & 529 & 659 \\
$k=4$ & 659&  659& 659 & 846& 678& 529& 846& 846& 846& 659\\
$k=5$ & 659&  659& 846 & 846& 678& 529 & 529& 529 & 846& 659 \\
$k=6$ & 659&  529& 659 & 846& 678& 846& 529& 846& 659& 846 \\
$k=7$ & 529&  529& 659 & 846& 659& 659& 529& 846& 659& 529\\
$k=8$ & 529& 678& 659& 846& 529& 846& 529 & 529& 846& 846 \\
$k=9$  & 529& 678& 659& 659& 529& 529& 529& 529& 659 & 846\\ \hline
\end{tabular}
\bigskip

{\sc Table}~2.1.~~Values of  total stopping time  $\sigma_\In (n)$ for
$n= n_0+10j +k,$
with $n_0 := 100\lfloor \pi \cdot 10^{35} \rfloor =
31, 415, 926,535, 897,932, 384,626,433, 832, 795, 028 ,800.$ \bigskip 
\end{center}


We observe that the total stopping time function takes only a few different values, 
namely: 529, 659, 678 and 846, and these four values occur intermixed in a 
somewhat random-appearing
way, but with some regularities. Note that around $n_0 \sim 3.14 \times 10^{37}$ the predicted
``average size" of a trajectory is $6.95212 \log n_0\approx 600$. In the data here we 
also observe  ``jumps"
of size between the occurring values on the order of $100$.

\smallskip
This is not a property of just this starting value.
In Table 2.2 we give similar data for blocks of $100$  near $n= 10^{35}$ and $10^{36}$,
respectively. Again we observe that there are also four or five values occurring,
but now they are different values.  
In this table we present data on  two other statistics: the  {\em frequency} statistic gives the 
count of these number of occurrences of each value,
and the {\em $1$-ratio} statistic denotes the fraction of odd iterates occurring in the given trajectory
up to and including when $1$ is reached. It is an experimental 
fact that all sequences in the table having  the same
total stopping time also 
have the same $1$-ratio.  In the first two blocks the value $\sigma_{\infty}(n)=481$ (resp. $351$)
that occurs with frequency $1$
 is that for the intial value $n =10^{35}$ (resp. $n=10^{36}$) in the given interval;
these initial values are  unusual in being divisible by a high power of $2$. 
Probabilistic models for the $3x+1$-function iteration  predict
that even and odd iterates will initially occur with equal frequency, so we may
anticipate the  $1$-ratio values to be relatively close to $0.5$.

\begin{center}
\begin{tabular}{|c|c|c||c|c|c||c|c|c|} \hline
~ & (a)~ $10^{35}$ & ~& ~&  (b)~ $10^{36}$ &~&~& (c)~$n_0$ & 
\\ \hline
 $ \sigma_\In (n)$ & freq. & $1$-ratio &  ~$ \sigma_\In (n)$ & freq. & $1$-ratio & $ \sigma_\In (n)$ 
 & freq.& $1$-ratio \\ \hline
481&  1  & 0.47817&351 & 1 &  0.41594 & 529 & 38& 0.48204\\
508 & 19& 0.48622 & 467& 72& 0.46895&659&28&0.51138 \\
573 & 49& 0.50261& 508 & 21& 0.48228&678&7&0.51474\\
592& 10&  0.50675 &519& 6  & 0.48554&846 & 27& 0.53782\\
836 & 21& 0.54306 &~& ~ &~ &&&\\
 \hline
\end{tabular}
\bigskip

{\sc Table}~2.2~~Values of total stopping time,  their frequencies, and $1$-ratio for \\
(a) $10^{35} \leq n \leq 10^{35}+99$, ~~~(b) $10^{36} \leq n \leq 10^{36}+99$,
(c)  $n_0 \leq n \leq n_0 +99$. 
\bigskip\\
\end{center}

The data in Table 2.2 suggests the following heuristic: as $n$ increases
only a few values of $\sigma_{\infty} (n)$ locally 
occur over short intervals; there is then a slow variation in which values of $\sigma_{\infty} (n)$ occur.
However these local values are separated
from each other by relatively large  ``jumps" in size. 
We stress that this is a purely empirical observation, nothing like this is 
rigorously proved! Our heuristic did not quantify what is a   ``short interval" 
and it did not quantify what ``relatively large jumps" should mean.
Even the existence of finite values for $\sigma_{\infty}(n)$ in the tables  presumes the
$3x+1$ conjecture is true for all numbers in the table.

 %
 %

\subsection{Probabilistic Models}

A challenging feature of the $3x+1$  problem is the huge gap between what can
be observed about its behavior in computer experiments and what can be rigorously proved.
Attempts to understand and predict features of empirical experimentation
have led to the following curious outcome: {\em the use of  probabilistic models to
describe a deterministic process! } This gives another theme of research
on this problem:
the construction and analysis of probabilistic and 
stochastic  models for various aspects of the iteration process. 

\smallskip
A basic probabilistic model of iterates of the $3x+1$ function $T(x)$ proposes
 that most trajectories for $3x+1$ iterates have equal numbers of
 even and odd iterates, and that the parity of
 successive iterates behave in some sense like independent
 coin flips. A key observation of Terras \cite{Ter76} and Everett \cite{Ev77}, leading
 to this model, is that the initial iterates of 
 the $3x+1$  function have this property (see Lagarias \cite[Lemma B]{Lag85}).) 
 This probabilistic model suggests that most trajectories plotted 
 on a logarithmic vertical scale should appear 
 close to  a  straight line having  negative slope 
 equal to  $-\frac{1}{2} \log \frac{3}{4} \sim -0.14384,$
 and should thus take about    $6.95212 \log n$ steps to reach $1$.
 
 \smallskip
 The corresponding behavior of iterates of the Collatz function $C(x)$ is more complicated.
 The allowed patterns of even and odd Collatz function iterates always 
 have an even iterate following  each odd iterate. 
 Probabilistic models taking this into account are more complicated to formulate and analyze
 than that for the $3x+1$ function; this is a main reason for studying the $3x+1$ function
 rather than the Collatz function. 
 Use of the probabilistic model above
 allows the heuristic inference that  Collatz iterates will be even about two-thirds of the time.

\begin{figure}\label{fig23c}
\centering
\includegraphics[width=3.5in]{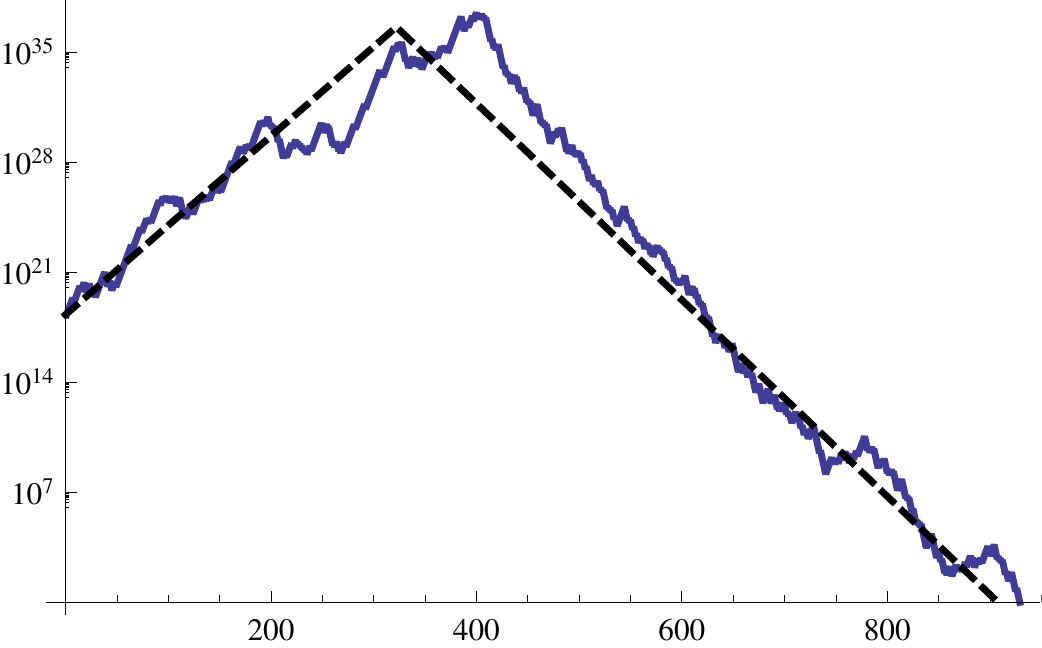}
\caption{
Extremal  trajectory
$n_1=1 980 976 057 694 878447$ 
given in Oliveira e Silva's Table 6.}
\end{figure}

\smallskip
A  variety of fairly complicated stochastic models,
many of which are rigorously analyzable (as probability models), have now been formulated to
model various aspects of these iterations, see Kontorovich and Lagarias \cite{KL09} (in this volume).
Rigorous results for such models lead to  
heuristic predictions for the statistical behavior of iterates of the generalized $3x+1$ map. 
The model above predicts the behavior of ``most" trajectories.
A small number of trajectories may exhibit quite different behavior.
 One may consider those trajectories that 
 that seem to offer maximal value of some iterate of $T^{(k)}(n)$ compared to $n$.
 Here a probabilistic model (see \cite[Sec. 4.3]{KL09} in this volume) predicts that the  statistic
 $$
 \rho(n)  := \frac{\log( \max_{k \ge 1} \left(T^{(k)}(n)\right)}{\log n}
 $$
 as $n \to \infty$ should have $\rho(n)  \le 2 + o(1)$ for all sufficiently large $n$. 
  Figure \ref{fig23c} offers a plot of the  trajectory, for 
  the value $n_1=1 980 976 057 694 878 447,$
which attains the largest value
of the statistic $\rho(n)$
 over $1 \le n \le 10^{18}$; 
this  value was found by Oliveira e Silva \cite[Table 6]{OeS09} (in this volume).
 This example has $\rho(n_1) \approx 2.04982$. Probabilistic models suggest that
 the extremal trajectories of this form will approach a characteristic shape
 which consists of two line segments, one of length $7.645 \log n$ steps of slope
 about $0.1308$ up to the maximal value
 of about $2 \log n$, the second of about $13.905 \log n$ steps of slope about $-0.1453$
 to 0, taking $21.55 \log n$ steps in all. This shape is indicated 
 by the dotted lines  on Figure \ref{fig23c} for comparison purposes.

%
%
%

Another  prediction of such stochastic models, relevant to the $3x+1$
conejcture, is that the  number of iterations required
for a positive integer $n$ to iterate to $1$ under the $3x+1$ function 
$T(x)$ is at most $41.677647 \log n$ (see \cite{LW92}, \cite[Sect. 4]{KL09}).
In particular such models predict, in a quantitative form, that
there will be no divergent trajectories.

\smallskip
These stochastic models can be generalized
to model the behavior of many generalized $3x+1$ functions, and they
make qualitatively different predictions depending on the function. For
example, such models  predict that no orbit of iteration of  the $3x+1$ 
function  ``escapes to infinity"
(divergent trajectory). 
However  for the {\em $5x+1$ function}  given by
$$
T_5(x) = 
\left\{
\begin{array}{cl}
\df{5x+1}{2} & \mbox{if}~ x \equiv 1~~ (\bmod ~2 ), \\
~~~ \\
\df{x}{2} & \mbox{if} ~~x \equiv 0~~ (\bmod~2),
\end{array}
\right.
$$
similar stochastic models predict that almost all orbits  should  ``escape to infinity"
(\cite[Sect. 8]{KL09}).
These predictions are  supported by experimental computer evidence, but   
it remains an unsolved problem to prove that there exists even one trajectory for the $5x+1$
problem that ``escapes to infinity".

 \smallskip
 There remains considerable research to be done on further developing
 stochastic models. 
 The experiments on the  $3x+1$ iteration reported above in \S3.2 exhibit some patterns
  not yet explained by stochastic models. 
  In particular,  the 
  behaviors of total stopping times observed   in Tables 2.1 and 2.2,
  and the heuristic presented there, 
  have not  yet been justified
  by suitable stochastic models.

%
%
%
\section{Generalized $3x+1$ functions}
\hsp

The original work on the $3x+1$ problem  viewed  it as a problem in number theory.
Much of the more  recent work views  it as an example of a  special kind of discrete dynamical
system, as exemplified by the lecture notes volume  of G. J. Wirsching \cite{Wir98}.
As far as generalizations are concerned, a very useful class of functions 
has proved to be  the 
set of generalized Collatz functions
which are defined below. These possess both number-theoretical 
and dynamic properties; the number-theoretic properties have to do
with the existence of $p$-adic extensions of these maps for various primes $p$. 

\smallskip
At present  the $3x+1$ problem is most often  viewed as a discrete dynamical
system of an arithmetical kind. It can then be treated as 
 a special case,  within the framework of a general class
of such dynamical systems. But what should be the correct degree of generality
in such a class?

\smallskip
There is significant interest in
exploring the behavior of dynamical systems of an arithmetic nature, 
since these may be viewed as ``toy models"
of more complicated dynamical systems arising
in mathematics and physics. 
There are a wide variety of interesting arithmetic dynamical systems. 
The book of  Silverman \cite{Sil07} studies
the iteration of 
algebraic maps on algebraic varieties. 
The book  of Schmidt \cite{Sch95}
considers dynamical systems of algebraic origin, meaning 
$\ZZ^d$-actions on compact metric groups, using ergodic theory and symbolic methods.
The book of Furstenberg \cite{Fur81} considers various
well structured arithmetical dynamical systems; for a further development 
see Glasner \cite{Gla03}.
The generalized $3x+1$ functions studied in this book 
provide  another  distinct type of  arithmetic discrete dynamical system. 

\smallskip
We present  a taxonomy of several classes of functions
which represent  successive generalizations of the $3x+1$ function. 
The simplest generalization of the $3x+1$ function
is  the $3x+k$ function,  which is defined for $k \equiv 1$ or $ 5 ~(\bmod 6)$, by 
$$
T_{3,k}(x) = 
\left\{
\begin{array}{cl}
\df{3x+k}{2} & \mbox{if}~ x \equiv 1~~ (\bmod ~2 ) ~, \\
~~~ \\
\df{x}{2} & \mbox{if} ~~x \equiv 0~~ (\bmod~2) ~.
\end{array}
\right.
$$
The generalization of the $3x+1$ conjecture to this situation is twofold: first, 
that under iteration every orbit becomes eventually periodic, and second,
that there are only a finite number of  cycles (periodic orbits). 
This class of  functions occurs in the study of cycles of the $3x+1$ function
(Lagarias \cite{Lag90}). Note that the $3x+1$
function $T(x)$ can be extended to  be well defined on the set of all rational numbers
having  odd denominator, and a rescaling of any $T$-orbit of such a rational number $r=\frac{n}{k}$
to clear its denominator $k$ will give an orbit of the map $T_{3, k}$.
Thus, integer cycles of the $3x+k$ function correspond to rational cycles of the
$3x+1$ function having denominator $k$.

\smallskip
 To further generalize, let  $d \ge 2$ be a fixed integer and consider the function
defined for integer inputs $x$ by
\beql{121}
f(x) = \frac{a_i x + b_i}{d} ~~~\mbox{if}~~x \equiv i~~(\bmod ~d), ~~0 \le i \le d-1,
\eeq
where  $\{ (a_i, b_i): 0 \le i \le d-1\}$ is a collection of integer
pairs. Such a function is called {\em admissible} 
if  the integer pairs $(a_i, b_i)$ satisfy the
condition
\beql{130}
i a_i + b_i \equiv 0 ~~(\bmod~d) ~~~\mbox{for}~~0 \le i \le d-1.
\eeq
This condition is  necessary and sufficient 
 for the map $f(x)$ to take integers to integers.
These functions $f(x)$ have been called 
{\em generalized Collatz functions,}  or 
 {\em $RCWA$ functions}
(Residue-Class-Wise  Affine functions). 
Generalized Collatz functions have the nice feature that they
have  a unique continuous extension to
the space $\ZZ_d$ of $d$-adic integers in the sense of Mahler \cite{Ma61}.

\smallskip
An important  subclass of generalized Collatz functions are
those of  {\em relatively prime type}.  These are the subclass of 
generalized Collatz
functions for which
 \beql{122} 
 \gcd ( a_0 a_1 \cdots a_{d-1}, d) =1. 
 \eeq
 This class includes the $3x+1$ function $T(x)$ but not the Collatz function $C(x)$ itself.
 It includes the $5x+1$ function $T_5(x)$, which as mentioned above 
 appears to have  quite different long-term dynamics  on the integers $\ZZ$
 than does the $3x+1$ function.
   Functions in this class have the additional property that their unique extension to the 
  $d$-adic integers $\ZZ_d$ has the $d$-adic Haar measure as an invariant measure.
  This permits ergodic theory methods to be applied to their study, see
 the survey paper of Matthews \cite[Thm. 6.2]{Mat09} (in this volume)
 for many examples.

\smallskip
 As a final generalization, one may consider the class of
integer-valued functions, which when restricted to 
residue classes $(\bmod \, d)$ are given by a polynomial $P_i(x)$
for each class $i~(\bmod \, d).$  
Members of this class of  functions have arisen in several places in mathematics.
They are now widely called {\em quasi-polynomial functions} or {\em quasi-polynomials}.
Quasi-polynomials  appear in commutative algebra and algebraic geometry, 
  in describing the Hilbert functions of certain semigroups, in
  a well known theorem of Serre, see
  Bruns and Herzog \cite[pp. 174--175]{BH93} and Bruns and Ichim \cite{BI07}.
In another direction, functions that  count the number of
lattice points inside dilated rational
polyhedra have been shown to be quasi-polynomial  functions (on the positive integers), 
starting with work of Ehrhart \cite{Ehr76}, see Beck and Robins \cite{BR07} and 
Barvinok \cite[Chap. 18]{Bar08}.
 They also have recently appeared in differential algebra in connection with
 $q$-holonomic sequences, see Garoufalidis \cite{Ga10}.
Such functions were introduced in group theory by 
G. Higman in 1960 \cite{Hi60} under the name  PORC functions  
(polynomial on residue class functions).
Higman's motivating problem was the enumeration of $p$-groups,
cf. Evseev \cite{Ev08}.
The class of all quasi-polynomial functions is closed under addition and pointwise
multiplication, and forms a commutative ring under these operations. 

\smallskip
We arrive at the following  taxonomy of function classes of increasing generality:
\begin{eqnarray*}
\{3x+1 ~\mbox{function} ~ T(x) \} & \subset &  \{3x+k ~\mbox{functions}~ T_{3,k}(x) \} \\
& \subset& \{ \mbox{generalized Collatz ~functions~of~relatively~prime~type} \} \\
&\subset & \{ \mbox{generalized Collatz~functions}\} \\
&\subset & \{ \mbox{quasi-polynomial ~functions}\} .
\end{eqnarray*}

For applications in mathematical logic, it
has proved useful to further widen the definition of generalized Collatz functions
to allow {\em partially defined functions}.  Such functions are 
 obtained by dropping the 
admissibility condition \eqn{130}; they map integers to
rational numbers having denominator dividing $d$. If a non-integer
value is encountered, then  one cannot iterate such a 
function further.
In this circumstance we  adopt the
convention that if a non-integer  iteration value is encountered, the calculation
stops in a special ``undefined"   state. This framework allows the encoding of 
partially-defined (recursive)
functions.  One can use this convention to also define composition of 
partially defined functions.

%
%
%
\section{Research Areas}
\hsp

Work on the $3x+1$ problem cuts across many  fields of mathematics.
Six  basic areas of research on the problem are:
(1) {\em number theory}: analysis of  periodic orbits of the map;  (2)  {\em dynamical
systems}: behavior of generalizations of the $3x+1$ map;  (3){\em  ergodic theory}:
invariant measures for generalized maps;  (4) {\em theory of computation}: undecidable
iteration problems;  (5) {\em stochastic processes and probability theory}:
models yielding heuristic predictions for the behavior of iterates; and 
(6) {\em  computer science}: algorithms
for computing iterates and statistics, and explicit computations.
We  treat these in turn.\\

(1) {\em Number Theory} \\

The connection with number theory is immediate: the $3x+1$ problem
 is a  problem in arithmetic,
whence it  belongs to elementary number theory. Indeed it is classified as an
unsolved problem in number theory by R. K. Guy \cite[Problem E16] {Guy04}.
The study of  cycles of the $3x+1$ map
leads to problems involving exponential Diophantine equations.
The powerful work of Baker and Masser--W\"{u}stholz on linear forms in logarithms  gives
information on the non-existence of cycles of various lengths  having specified
patterns of even  and odd iterates. 
A class of generalized $3x+1$ functions has been defined in a number theory framework, 
in which arithmetic operations on  the domain of  integers are replaced with
such operations on  the ring of integers of
an algebraic number field, or by  function field analogues
such as a polynomial ring with coefficients in a finite field.  Number-theoretic results
are surveyed in the papers of Lagarias \cite{Lag85} and Chamberland \cite{Cha07}
 in this volume. \\

(2) {\em Dynamical Systems} \\

The theory of discrete dynamical systems concern
the behavior  of functions under iteration; that of
continuous dynamical systems concern flows or solutions to 
differential equations.
 The $3x+1$ problem 
can be viewed as 
iterating a map, therefore it is  a discrete dynamical system on
the state space $\ZZ$. This viewpoint was taken in
Wirsching \cite{Wir98}.
The important operation for iteration is {\em composition of functions}. 
One  can formulate iteration and composition questions in the
general context of universal algebra, 
cf. Lausch and Nobauer \cite[Chap. 4.5]{LN73}.
In the taxonomy above, the classes 
of  generalized $3x+1$ functions,  and quasi-polynomial 
functions are each  closed under addition and composition of functions. 
The iteration properties of  the first three classes of functions above have been studied,
in connection with the $3x+1$ problem and the theory
of computation. However the iteration of general quasi-polynomial  functions
remains an unexplored research area. 

\smallskip
Viewing the problem this way suggests that it would be useful
in the study  of the $3x+1$ function to 
obtain dynamical systems on larger domains,
including the real numbers $\RR$ and the complex numbers $\CC$.
Other extensions include defining analogous functions on  the ring $\ZZ_2$ of $2$-adic integers,
or,  for generalized $3x+1$ maps,  on a ring of $d$-adic integers, for a value of $d$
determined by the function.  When one considers generalized 
$3x+1$ functions on larger domains, a wide variety of 
behaviors can occur. These topics are considered in the papers of Chamberland \cite{Cha07}
and Matthews \cite{Mat09} in this volume. For a general framework on topological
dynamics see Akin \cite{Ak93}. \\

(3) {\em Ergodic Theory} \\

The connection with ergodic theory arises as an outgrowth of
the  dynamical systems
viewpoint, but adds the requirement of the presence of
an invariant measure. It was early observed
that there are finitely additive  measures which
are preserved by the $3x+1$ map on the integers.
Extensions of  generalized $3x+1$ functions  to $d$-adic integers lead to maps
invariant under standard measures (countably additive measures). 
For example, the  (unique continuous) extension  of the $3x+1$  map to the $2$-adic integers
has $2$-adic measure as an invariant measure, and the map is ergodic
with respect to this measure.  
 Ergodic theory  topics are considered in the surveys of 
Matthews  \cite{Mat09} and Kontorovich and Lagarias \cite{KL09} in this volume. 
An interesting open problem is to classify all invariant measures for generalized
$3x+1$ functions on the $d$-adic integers. \\

(4) {\em Mathematical Logic and the Theory of Computation} \\

The connection to logic and the theory of computation starts with the result of 
Conway that there is a generalized $3x+1$ function whose iteration can
simulate a universal computer. Conway \cite{Con72} exhibited an unsolvable iteration
problem for a particular generalized $3x+1$  function: starting with a given input which is a 
positive integer $n$, decide whether or not some iterate of this map with
this input is ever a power of $2$.  
In this connection note  that the $3x+1$ problem can be reformulated
as asserting that, starting from any positive integer $n$, some iterate $C^{(k)}(n)$ of
the Collatz function (or of the $3x+1$ function) is a power of $2$.
It turns out that iteration of $3x+1$-like
functions had already been considered in understanding the power of 
some logical theories even
in the late 1960's; these involved partially defined functions taking integers
to integers (with undefined output for some integers), cf. Isard and Zwicky \cite{IZ70}.
 More recently
such functions have arisen in studying the computational power of
``small'' Turing machines, that are too small to encode a universal computer.
These topics are surveyed in the paper of Michel and Margenstern \cite{MM09} in this volume.\\

(5) {\em Probability Theory and Stochastic Processes}\\

A connection to probability theory and stochastic processes arises when
one attempts to model the behavior of the $3x+1$ iteration on large sets of
integers. This leads to heuristic probabilistic models for the iteration,
which allow predictions of its behavior. Some authors have argued that
the iteration can be viewed as a kind of pseudo-random number generator,
viewing the input as being given by a probability distribution, and then
asking how this probability distribution evolves under iteration.
In the reverse direction, one can study trees of inverse iterates (the inverse
map is many-to-one, giving rise to a unary-binary tree of inverse iterates).
Here one can ask for facts about the structure of such trees whose root node 
is an integer picked from some probability distribution. One can
model this by a stochastic model corresponding to random tree growth, e.g. a branching
random walk. These topics 
are  surveyed in the paper of Kontorovich and Lagarias \cite{KL09} in this volume. \\

(6) {\em Computer Science: Machine Models, Parallel and Distributed Computation} \\

In 1987 Conway \cite{Con87} (in this volume) formalized the Fractran model of 
computation as a universal computer model,
based on his earlier work related to the $3x+1$ problem. 
This computational model is related to
the register machine (or counter machine) model of Marvin Minsky (\cite{Min61},
\cite[Sect. 11.1] {Min67}). Both these  machine models have recently been seen as
relevant for developing models of
computation using chemical reaction networks, and to biological computation,
see   Soloveichik et al \cite{SCWB08}  and Cook et al. \cite{CSWB10}.

The necessity to make computer experiments to test the $3x+1$ conjecture,
and to explore various properties and patterns of the $3x+1$ iteration, leads
to other questions in  computation.
One has the research problem of developing efficient 
algorithms for computing on a large scale, using either
 parallel computers or a distributed computer system. 
 The $3x+1$ conjecture  has been tested to a very large value of $n$,
see the paper of Oliveira e Silva \cite{OeS09} in this volume.
The computational method 
 used in \cite{OeS09}  to obtain record results
 can be parallelized. 
 Various large scale computations for the $3x+1$ problem have used
distributed computing, cf. Roosendaal \cite{Roo}.  \\



%
%
%
 
\section{Current Status}

We  give a brief summary  of the current status of the problem,
 which further elaborates 
answers to the two  questions raised in the introduction. 
%
%
%

\subsection{Where does research currently stand on the $3x+1$ problem?}

The $3x+1$ problem remains unsolved, and a solution remains unapproachable
at present.  
To quote a still valid dictum of Paul Erd\H{o}s (\cite[p. 3]{Lag85}) on the problem:   

\begin{quotation}
``Mathematics is not yet ready for such problems."
\end{quotation}

\smallskip
Research has established  various  ``world records", all of which rely on 
large computer calculations (together with various theoretical developments). 

\begin{enumerate}
\item[(W1)]
The $3x+1$ conjecture has now been verified for all 
$n < 20 \times 2^{58} \approx 5.7646 \times 10^{18}$
(Oliveira e Silva \cite{OeS09} (in this volume)).
\item[(W2)]
The trivial cycle $\{1, 2\}$ is the only cycle of the $3x+1$ function on the positive
integers having 
period length less than $10,439,860,591$. It is also the only 
cycle containing less than $6,586,818,670$ odd integers 
(Eliahou \cite[Theorem 3.2]{Eli93}\footnote{This number is the bound $(21,0)$ 
given in  \cite[Table 2]{Eli93}. The smaller values in Table 2 are 
now ruled out by the computations in item (W1) above.}).

\item[(W3)]
Infinitely many positive integers $n$ take at least $6.143 \log n$ steps to reach $1$
under iteration of the $3x+1$ function $T(x)$ (Applegate and Lagarias \cite{AL03}).

\item[(W4)]
The positive integer $n$ with the largest currently known value of $C$, such that it takes $C \log n$
iterations of the $3x+1$ function $T(x)$ to reach $1$, is 
$n=7,219,136,416,377,236,271,195$
with $C\approx 36.7169$ (Roosendaal \cite[$3x+1$ Completeness and Gamma records]{Roo}). 

\item[(W5)]
The number of integers $1 \le n \le X$ that iterate to $1$ is at least $X^{0.84}$,
for all sufficiently large $X$ (Krasikov and Lagarias \cite{KL03}).
\end{enumerate}

There has also been considerable progress made on showing the nonexistence of various
kinds of periodic points for the $3x+1$ function, see Brox \cite{Br00} and 
Simons and de Weger \cite{SdW05}.
 These bounds are based on number-theoretic
methods involving Diophantine approximation. \\

%
%
%

\subsection{ Where does research  stand on generalizations of the $3x+1$ problem?}

It has proved fruitful to view the $3x+1$ problem as a special case of 
wider classes of functions. These function classes appear naturally
as the correct level of generality for basic results on iteration; 
this resulted in the taxonomy of function classes given in \S3.
There are some general results for these classes and  many unsolved problems. 

\smallskip
The $3x+k$ problem seems to be the correct level of generality for studying rational cycles
of the $3x+1$ function (\cite{Lag90}). There are extensive results on cycles of the
$3x+1$ function, and the methods generally apply to the $3x+k$ function as
well, see the 
survey of Chamberland \cite{Cha07} (in this volume).

\smallskip
The class of  generalized $3x+1$  functions of relatively prime type is
a very natural class from the ergodic theory viewpoint, since this is the class on which
the  $d$-adic extension
of the function has $d$-adic Haar measure as an invariant measure.
The paper of Matthews \cite{Mat09} (in this volume) reports 
general ergodicity results and raises many questions about
such functions.

\smallskip
The class of generalized Collatz functions  has the property that
all functions in it  have a unique continuous extension to 
the domain of $d$-adic integers $\ZZ_d$.  
This general class is known to contain undecidable iteration problems, as 
discussed in the paper of Michel and Margenstern \cite{MM09} (in this volume). 
The dynamics of general functions in this class is only starting to be explored; many 
 interesting examples are given in the paper of Matthews \cite{Mat09} (in this volume).
An interesting area worthy of  future development is  that of determining the existence
and structure of invariant Borel measures for such functions on $\ZZ_d$, and 
determining whether there is some relation
of their structure to undecidability  of the associated iteration problem.

%
%
%

\subsection{How can this be a hard problem, when it is so easy to state?}

Our answer is that there are two different mechanisms yielding hard problems, either
or both of which 
may apply to the $3x+1$ problem. The first is ``pseudorandomness"; this involves
a connection with ergodic theory. 
The second is ``non-computability". Both of these are discussed in detail in this
volume. 

\smallskip
The ``ergodicity" connection has been independently noted by a number of  people, see for example
Lagarias \cite{Lag85} (in this volume) and Akin \cite{Ak04}. 
The unique continuous extension
of the $3x+1$ map $T(x)$  to the $2$-adic integers $\ZZ_2$ gives a function which is known to be
ergodic in a strong sense, with respect to the $2$-adic measure. It is topologically
and metrically conjugate
to the shift map, which is a maximum entropy map. 
The iterates of the shift function are completely unpredictable in the ergodic theory sense. Given
a random starting point, predicting the parity of the $n$-th iterate 
for any $n$ is a ``coin flip" random variable. 
The $3x+1$ problem
concerns the behavior of iterating this function on the set of integers $\ZZ$, which is
a dense subset of  $\ZZ_2$, having 
$2$-adic measure zero. The difficulty is then in finding and understanding 
non-random regularities in the iterates when restricted to $\ZZ$. Various probabilistic
models  are discussed in the paper of Kontorovich and Lagarias \cite{KL09} (in this volume). 
Empirical evidence
seems to indicate that the $3x+1$ function on the domain $\ZZ$ 
retains the ``pseudorandomness"
property on its initial iterates until the iterates enter a periodic orbit. 
This supports the $3x+1$ conjecture
and at the same time deprives us of any obvious mechanism to prove it, since mathematical arguments
exploit the existence of structure, rather than its absence. 


\smallskip
A connection of a generalized Collatz function
to ``non-computability" was  made by Conway \cite{Con72} (in this volume),
as already mentioned. 
Conway's undecidability
 result indicates that the $3x+1$ problem could be close to the unsolvability threshold.
It is currently unknown whether the $3x+1$ problem is itself undecidable, however no method is
currently known to approach this question. 
The survey of Michel and Margenstern \cite{MM09}  (in this volume) describes many results on
generalized $3x+1$ functions that exhibit  undecidable or difficult-to-decide 
iteration problems.
The $3x+1$ function might conceivably belong to a smaller class of
generalized $3x+1$ functions that evade  undecidability results that encode
universal computers. Even so, it conceivably might  encode an undecidable problem, arising
by another (unknown) mechanism. 
As an example, could the following question be undecidable:
``Is there any positive integer $n$ such that
 $T^{(k)}(n) >1$ for $1 \le k \le 100 \log n$?"

%
%
%

\section{Hardness of  the $3x+1$ problem}\label{sec7}

\smallskip
Our viewpoint on hard problems has evolved since 1900,
starting  with Hilbert's
program in logic and proof theory and 
benefiting from  developments in the theory of computation.
 Starting in the 1920's, Emil Post uncovered
great complexity in studying some very simple computational problems, 
now called ``Post Tag Systems".
 A {\em Tag system} in the class $TS(\mu, \nu)$ consists of a set of rules 
for transforming words  using letters from an alphabet $\sA = \{ a_1, ..., a_{\mu} \} $ of $\mu$ symbols, 
a  deletion number (or shift number) $\nu \ge 1$, 
and  a set of $\mu$ production rules 
$$
a_j \mapsto w_j:= a_{j,0} a_{j,1} \cdots a_{j, n_j},~~~ 1 \le j \le \mu,
$$
 in which the output $w_j$ is a finite string (or  word) of length $n_j$ in the alphabet $\sA$. 
Starting from an initial string $S$ a Tag system looks at the leftmost symbol of $S$, call it $a_j$, 
then attaches to the right end of the string the word $w_j$, 
and finally deletes the first $\nu$ symbols
of the resulting string $S w_j$, thus obtaining a new string $S'$. Here the ``tag" is the set of
symbols $w_j$ attached to the end of the word, and the iteration halts if a word of length
less than $\nu$ is encountered. 
The {\em halting problem} is the question
of deciding whether for an arbitrary initial word $S$, iteration 
eventually reaches the empty word.
The {\em reachability problem}
is that of deciding whether, given words $S$ and $\tilde{S}$, starting from word $S$ will
ever produce word $\tilde{S}$ under iteration. The halting problem is a special
case of the reachability problem.
Post \cite{Pos65} reports that in 1920--1921 he found a complete decision
procedure\footnote{Post did not  publish his proof. A decision procedure for both problems is
outlined in de Mol \cite{DM07}.}
for 
the case $\mu=2, \nu=2$, i.e. the class $T(2,2)$. 
He then tried to solve the case $\mu=2, \nu>2$,
 without success. He reported \cite[p. 372]{Pos65}
  that the special case $\mu=2, \nu=3$  with $\sA =\{0, 1\}$
 and the two production rules
 \beql{730a}
 0 \mapsto w_0=00, ~~~ 1 \mapsto w_1=1101
 \eeq
already  seemed to be an intractable problem. We shall term this problem

\medskip
{\bf Post's Original Tag Problem.} {\em Is there a recursive decision procedure for
the halting problem for the Tag system  in $T(2, 3)$ given by the rules $0 \mapsto 00$
and $1 \mapsto 1101$?}

\smallskip
\noindent Leaving this question aside, Post considered the parameter range $\mu>2, \nu=2$. He wrote
 \cite[p. 373]{Pos65}: 

\begin{quotation} 
For a while the case $\nu=2, \, \mu>2$
  seemed to be more promising, since it seemed to offer a greater chance of
 a finitely graded series of problems. But when this possibility was explored in
 the early summer of 1921, it rather led to an overwhelming confusion of classes of
 cases, with the solution of the corresponding problem depending more and more
 on problems of ordinary number theory. Since it had been our hope that the known
 difficulties of number theory would, as it were, be dissolved in the particularities of
 this more primitive form of mathematics, the solution of the general problem of
 ``tag" appeared hopeless, and with it our entire program of the solution of 
 finiteness problems.
  \end{quotation}

\smallskip
\noindent  Discouraged by this, Post reversed course and went on to
obtain  a ``Normal Form Theorem" (\cite{Pos43}), published  in the 1940's,  showing
that a general logical problem could be reduced to a form
slightly more complicated than Tag Systems. 
In 1961 Marvin
Minsky \cite{Min61}   proved that Post Tag Systems were undecidable problems in general.
In the next few years  Hao Wang \cite{Wang63},  J. Cocke and M. Minsky \cite{Coc64} and S. Ju. Maslov \cite{Mas64}
independently showed undecidability for the subclass of Post Tag Systems
consisting of those with $\nu=2$, thus showing that Post was right to quit trying to
solve problems in that class.  At present the recursive solvability or unsolvability in the class 
$T(2, \nu)$ remains open for all $\nu>2$. Post's original tag problem, which is 
the halting problem for one special function in $T(2,3)$,  is still unsolved,
see Lisbeth De Mol \cite{DM06}, \cite[p. 93]{DM08}, and for further work \cite{DM07}, \cite{DM09}.

\smallskip
Recently    de Mol showed that the $3x+1$ problem can
be encoded as a reachability problem for a tag system  in $T(3,2)$  (\cite[Theorem 2.1]{DM08}).
This tag system encodes the $3x+1$ function, and the reachability problem is:

\smallskip
{\bf $3x+1$ Tag Problem.} {\em
Consider the tag system $T_C$ in $T(3,2)$ with alphabet  $\sA=\{ 0, 1, 2\}$,
deletion number $\nu=2$,  and
production rules
$$
0\mapsto 12, \, 1 \mapsto 0, \, 2 \mapsto 000.
$$
 For each $n \ge 1$, if one starts from the configuration $S= 0^n$,  will
the tag system iteration for $T_C$ always  reach state $\tilde{S}=0$?}

\smallskip
In 1931 Kurt G\"{o}del \cite{God31} showed the existence of undecidable
problems:  he showed that  certain propositions were undecidable
in any logical system complicated enough to include elementary number
theory. This result showed that Hilbert's proof theory program could not be carried
out. Developments in the theory of computation showed that one of G\"{o}del's incompleteness
results corresponded to the unsolvability of the halting problem for Turing machines. 
This was based on the existence of a universal Turing machine, that could
simulate any computation, and in his 1937 foundational paper Alan Turing \cite{Tur37}
already showed one could be constructed of a not very large size. 

\smallskip
We now have a deeper appreciation of exactly how simple a problem can be and
still simulate a universal computer. Amazingly simple  problems of this sort have been found 
in recent years. Some of these involve cellular automata, a model of computation
developed by John von Neumann and Stansilaw M. Ulam in the 1950's. 
One of these  problems concerns  the possible behavior of  a very  simple one-dimensional
nearest neighbor cellular automaton,  Rule 110, using
a nomenclature introduced by Wolfram \cite{Wol83}, \cite{Wol84}. 
This rule was conjectured by Wolfram
 to give a universal computer (\cite[Table 15]{Wol87}, \cite[pp. 575--577]{Wol94}).
 It was proved  to be weakly universal by M. Cook (see  Cook \cite{Co04},
 \cite{Co09}). Here weakly universal means that the initial configuration of the
 cellular automaton is required to be ultimately periodic, rather than finite.
Another is John H. Conway's game of ``Life,"  first announced in 1970
in Martin Gardner's column in Scientific American  (Gardner \cite{Gar70}),
which is a two-dimensional cellular automaton,
having nearest neighbor interaction rules of a particularly simple nature.
Its universality as a computer was later established, see 
Berkelamp, Conway and Guy \cite[Chap. 25]{BCG04}. Further remarks on the
size of universal computers are given in the survey of Michel and Margenstern  \cite{MM09}
(in this volume). 

\smallskip
There are, however,  reasons to suspect that the $3x+1$ function is not
 complicated enough to be universal, i.e. to allow the encoding
of a universal computer in its input space. 
First of all, it is so simple to state that there seems very little room in it to encode the
elementary operations needed to create a universal computer. Second,
the $3x+1$ conjecture asserts that the iteration halts on the domain of
all positive integer inputs, so for each integer $n$, the value $F(n)$ of the
largest integer observed before visiting $1$ is recursive.
To encode a universal computer, one needs to represent all recursive functions,
including functions that grow far faster than any given recursive function $F(n)$.
It is hard to imagine how one can encode it here  as a question
about the iteration, without enlarging the domain of inputs.  Third, the $3x+1$ function
possesses the  feature that there is a nice (finitely additive) invariant measure
on the integers, with respect to which it  is completely mixing under iteration.
This is the measure that assigns mass $\frac{1}{2^n}$ to each complete arithmetic
progression $(\bmod~ 2^n)$, for each $n \ge 1$. This  fundamental observation
was made in 1976 by Terras \cite{Ter76}, and independently by Everett \cite{Ev77}
in 1977,
see Lagarias \cite[Theorem B]{Lag85} for a precise statement.
This  ``mixing property" seems to fight against the amount of organization needed
to encode a universal computer in the inputs.  We should caution that this 
observation by itself does not 
rule out the possibility that, despite this mixing property,
 a universal computer could  be encoded in a very
thin set of input values (of ``measure zero"), compatible with an invariant
measure.  It just makes it seem
difficult to do.  Indeed, the 1972 encoding of
a universal computer in the iteration of a certain
generalized $3x+1$  function
found by Conway \cite{Con72} (in this volume)
has the undecidability encoded in the iteration of a very thin set of integers.
However Conway's framework is different from the $3x+1$ problem
in that the halting function he considers is partially defined.


\smallskip 
Even if iteration of the $3x+1$ function
is not universal, it could still potentially be unsolvable. 
Abstractly, there may exist in an  axiomatic system 
 statements $F(n)$ for a positive integer predicate, such that
$F(1), F(2), F(3), ...$ are provable in the system for all integer $n$,  but the
statement $(\forall n) \,F(n)$ is not provable within the system. 
For example, one can let $F(n)$ encode a statement that there
is no contradiction in a system obtainable by a proof of length
at most $n$.  If the system is consistent,
then $F(1), F(2), ...$ will all individually be provable.
The statement $(\forall n)\, F(n)$ then encodes the
consistency of the system. But the consistency of a system sufficiently
complicated to include elementary number theory cannot be
proved within the system, according to   G\"{o}del's
second incompleteness theorem. 


\smallskip
The  pseudo-randomness or ``mixing" behavior of the $3x+1$ function also
seems to make it extremely resistant to analysis. 
If one could rigorously show a sufficient amount of mixing
is guaranteed to occur, in a controlled number of iterations in terms of the input size $n$,
then one could settle part of the $3x+1$ conjecture, namely prove the non-existence
of divergent trajectories. Here we have the fundamental difficulty of proving
in effect that the iterations actually do have an explicit  pseudo-random property. 
Besides this difficulty, there remains a second fundamental
difficulty: solving the number-theoretic problem 
 of  ruling out the existence of an enormously long non-trivial 
cycle of the $3x+1$ function. This problem
also  seems unapproachable  at present by known methods
of number theory. However the finite cycles problem does
admit proof of partial results, showing  the nonexistence
of non-trivial cycles having particular patterns of even and odd iterates. 

\smallskip
A currently  active and  important general area of research concerns the construction of
pseudo-random number generators:  these are deterministic recipes that
produce apparently random outputs (see Knuth \cite[Chap. 3]{Kn81}). More precisely, one is interested
in methods that take as input $n$ truly random bits and deterministically produce as output
$n+1$ ``random-looking" bits. These bits are to be ``random-looking" in the sense
that they appear random with respect to a given family of statistical tests,
and the output is then said to be pseudo-random with respect to this
family of tests. Deciding whether pseudo-random number generators exist
for statistical tests in various complexity classes is now seen as a fundamental question 
in computer science, related to the $P=NP$ probem, see
for example Goldreich \cite{Gol01}, \cite{Gol10}.
It may be that  resolving  the issue of the pseudo-random character of iterating the $3x+1$
problem will require shedding light on the general existence problem for
pseudo-random number generators.

\smallskip
All we can say at present is that the $3x+1$ problem appears very hard indeed. It  now seems
less surprising than it might have once seemed 
 that a problem as simple-looking as this one could be genuinely
difficult, and inaccessible to known methods of attack.

%
%
%
%

\section{Future Prospects}


We  observe first that further
improvements are surely possible on the ``world records"  (W1)--(W5) above. 
 In particular, concerning (W3), it seems  
 scandalous that  it is not known whether or not there are infinitely
many positive integers $n$ which iterate to $1$ under the $3x+1$ map $T(x)$
and  take at least the ``average" number $\frac{2}{\log 4/3} \log n \approx 6.95212 \log n $ steps 
to do so. Here the stochastic models for the $3x+1$ iteration
predict that at least half of all positive integers should
have this property! These ``world records" are  particularly worth improving if they can
shed more light on the problem. 
This could  be the case for 
world record (W5), where there is an underlying structure for obtaining lower bounds on
the exponent, which involves an 
 infinite family of nonlinear programs of increasing complexity (\cite{KL03}). 

\smallskip
Analysis of the $3x+1$ problem has resulted in the formulation of a large set of
``easier'' problems.  At first glance some of these 
 seem  approachable, but they also  remain unsolved,
 and are apparently difficult.  As samples, these include:

\begin{enumerate}
\item[(C1)]
({\em Finite Cycles Conjecture})
Does the  $3x+1$ function have 
 finitely many cycles
(i.e. finitely many  purely periodic orbits on the integers)?
This is conjectured to be the case.

\item[(C2)]
({\em Divergent Trajectories Conjecture-1})
Does the $3x+1$ function have a divergent trajectory,
i.e., an integer starting value whose iterates are unbounded?
This is conjectured {\em not} to be the case.

\item[(C3)]
({\em Divergent Trajectories Conjecture-2})
Does the $5x+1$ function have a divergent trajectory? This
is conjectured to be the case. 

\item[(C4)]
({\em Infinite Permutations-Periodic Orbits Conjecture})
If a  generalized Collatz function permutes the integers and 
is not globally of finite order, is it true that it has only finitely many periodic orbits?
The original Collatz function $U(n)$, which is a permutation,
 was long ago conjectured to have finitely many cycles.  A  conjecture 
  of this kind,  imposing extra conditions  on the permutation, 
 was formulated by Venturini  \cite[p. 303 top]{Ven97}.

\item[(C5)]
({\em Infinite Permutations-Zero Density Conjecture})
If a  generalized Collatz function permutes the integers, is it true that
every orbit has a (natural) density? 
Under some extra hypotheses 
 one may conjecture
that  all such orbits  have density zero;  compare Venturini \cite[Sec. 6]{Ven97}.
 \end{enumerate}

 \smallskip
Besides these conjectures, there also exist  open problems which may be more
accessible.
 One of the most intriguing of them concerns establishing 
lower bounds
for the number $\pi_1(x)$ of integers less than $x$ that get to $1$ under   
the $3x+1$ iteration. As mentioned earlier it is known
 (\cite{KL03}) that there is a positive constant $c_0$
such that 
$$
\pi_1(x) > c_0 x^{0.84}.
$$
It remains an open problem to show  that for each $\epsilon>0$ there exists a positive constant $c(\epsilon)$
such that
$$
\pi_1(x) > c(\epsilon) x^{1- \epsilon}.
$$
Many other  specific, but difficult,  conjectures for study can be found in the papers in this volume,
starting with the problems listed in Guy \cite{Guy83}.

\smallskip
 We now raise some  further research directions, related to the papers
 in this volume.  A first research direction  is to extend the class of functions for which the Markov models 
 of Matthews \cite{Mat09} can
 be analyzed. Matthews shows that the class of generalized $3x+1$ functions of relatively
 prime type (\cite[Sec. 2]{Mat09}) is  analyzable.
 He  formulates some conjectures for exploration.
  It would be interesting  to  characterize the possible $d$-adic invariant measures for
 arbitrary generalized Collatz functions.  It may be necessary to restrict to subclasses of
 such functions in  order to obtain nice characterizations. 

\smallskip
 A second research direction concerns the class of generalized $3x+1$ functions whose
 iterations extended to the set of $d$-adic integers are ergodic with
 respect to the $d$-adic measure, cf. Matthews \cite[Sec. 6]{Mat09}).
 
 \smallskip
{\bf Research Problem.} {\em  Does the class of generalized 
Collatz functions of relatively prime type contain a
 function which is ergodic with respect to the standard $d$-adic measure, whose iterations
 can simulate a universal computer? 
 Specificially, could it have  an unsolvable iteration problem of the form: 
  ``Given positive integers $(n, m)$ as input, does there
 exist $k$ such that the $k$-th iterate $T^{(k)}(n)$ equals $m$?"
Or does  ergodicity of the iteration preclude the possibility of simulating
universal computation?}

 \smallskip
  A third research direction concerns the fact
  that generalized Collatz 
   functions have now been found in many 
 other mathematical structures, especially if one generalizes further to integer-valued
 functions that are piecewise polynomial on residue classes $(\bmod~d)$. 
 These functions are the  quasi-polynomial  functions noted above, and they show up 
 in a number of algebraic contexts, particularly in counting lattice points in various regions.
 It may prove worthwhile to study the iteration of various special classes of 
 quasi-polynomial  functions arising in these algebraic contexts.
 
 

\smallskip
 At this point in time, in view of the intractability of problems (C1)--(C5) 
  it  also seems  a sensible task
 to formulate a new  collection of even simpler ``toy problems",
 which may potentially be approachable.  These may involve 
 either changing the problem or
 importing it into new contexts.
 For example, there appear to be accessible open problems concerning variants of the problem
 acting on finite rings (Hicks et al. \cite{HMYZ08}).
Another  promising recent direction  is the connection
 of these problems with generating sets for multiplicative arithmetical
 semigroups, noted by Farkas \cite{Far05}. 
 This has led to a family of more accessible problems, where various results can
 be rigorously established (\cite{AL06}). Here significant unsolved problems remain
 concerning the structure of such arithmetical semigroups.
 Finally it may prove profitable to continue
 the study, initiated by Klarner and Rado \cite{KR74}, of sets of integers
 (or integer vectors) closed under the action of a finitely generated semigroup of affine maps.

%
%
%
\section{Is the $3x+1$ problem a ``good" problem?}

There has been much discussion of what constitutes a good mathematical problem.
We cannot do  better than to recall the discussion of Hilbert \cite{Hi00} 
in his famous 1900 problem list. On the importance of problems
he said (\cite[p. 437]{Hi00}):

\begin{quotation}
The deep significance of certain problems for the advance of mathematical science
in general, and the important role they play in the work of the individual investigator,
are not to be denied. As long as a branch of science offers an abundance of problems,
so long is it alive; a lack of problems foreshadows extinction or the cessation of 
independent development. Just as every human undertaking pursues certain
objects, so also mathematical research requires its problems. It is also by
the solution of problems that the investigator tests the temper of his steel; he finds new
methods and new outlooks, and gains a wider and freer horizon.
\end{quotation}

\smallskip
Hilbert puts forward three  criteria that a  good mathematical problem ought to  satisfy: 

\begin{quotation}
It is difficult and often impossible to judge the value of a problem correctly in advance; for
the final award depends upon the gain which science obtains from the problem. Nevertheless
we can ask whether there are general criteria which mark a good mathematical problem.
An old French mathematician said: ``A mathematical theory is not to be considered
complete until you have made it so clear that you can explain it to the first man that you
meet on the street." This clearness and ease of comprehension, here insisted on for
a mathematical theory, I should still more demand for a mathematical problem if it is
to be  perfect; for what is clear and easily comprehended attracts, the complicated
repels us.

Moreover a mathematical problem should be difficult in order to entice us, but not
completely inaccessible, lest it mock at our efforts. It should be to us a guide post on
the mazy paths to hidden truths, and ultimately a reminder of our pleasure in its
successful solution.
\end{quotation}

\smallskip
 From the viewpoint of the Hilbert criteria for a good problem, we see that:

\medskip
(1)  The $3x+1$ problem is a clear, simply
stated problem;   

\medskip
(2) The $3x+1$ problem  is a difficult problem; 
 
 \medskip
 (3) The $3x+1$ problem  initially seems accessible, in that it  possesses a fairly intricate
 internal structure.  

\smallskip
But -- and it is a big ``but" -- 
the evidence so far suggests 
that obtaining  a proof of the $3x+1$ problem is   inaccessible!
Not only does this goal appear inaccessible, but various simplified conjectures derived
from it  appear to be completely inaccessible in their turn, leading to a regress to 
formulation of a series of 
simpler and simpler inaccessible problems, namely conjectures (C1)--(C5) listed in \S8.

\smallskip
 We conclude that 
 the $3x+1$ problem  comes close  to being a  ``perfect" problem in the Hilbert sense. 
 However it  seems to  fail   the last of Hilbert's requirements:
  It mocks our efforts! It is possible to work hard
 on this problem to  no result. It is definitely a dangerous problem! 
 It could well be that the $3x+1$ problem remains out of human reach. But maybe not.
Who knows?  \\

%

%
%
%
\section{Working on   the $3x+1$ problem}

Whether or not the $3x+1$ problem is a ``good" problem, it
 is not going away, due to its extreme accessibility. 
It offers a large and tantalizing 
variety of  patterns  in computer experiments.  This problem 
stands as a mathematical challenge for the 21-st century.

\smallskip
In working on this problem, the  most cautious advice, following Richard Guy \cite{Guy83} 
is: 

\begin{quotation}
 {\em  
 Don't try to solve these problems!
 }
 \end{quotation}
 
 \noindent But, as Guy said \cite[p. 35]{Guy83}, 
 some of you may be already scribbling, in spite of the warning! \medskip

 We also note  that Paul
Erd\H{o}s  said, in conversation, about its difficulty (\cite{Erd90}):

\smallskip
\begin{quotation}
{\em ``Hopeless. Absolutely hopeless."}
\end{quotation}

\smallskip
\noindent In Erd\H{o}s-speak, this means that
there are no known methods of approach
which gave any promise of solving the problem.
For other examples of Erd\H{o}s's use of the term ``hopeless" 
see Erd\"{o}s and Graham \cite[pp. 1, 27, 66, 105]{EG80}.

\smallskip
At this point we  may  recall  further advice of David Hilbert  \cite[p. \, 442]{Hi00} about  problem solving:

\begin{quotation}
If we do not succeed in solving a mathematical problem, the reason frequently
consists in our failure to recognize the more general standpoint from which the
problem before us appears only as a single link in a chain of related problems.
After finding this standpoint, not only is this problem frequently more accessible
to our investigation, but at the same time we come into possession of a method
that is applicable to related problems.
\end{quotation}

\smallskip
\noindent The quest for 
generalization cuts in two directions,
for Hilbert also says \cite[p. \,442] {Hi00}:

\begin{quotation}
He who seeks for methods without having  a definite problem in mind seeks
for the most part in vain.
\end{quotation}

\smallskip
Taking this advice 
into account,  researchers  have treated  many generalizations of the $3x+1$ problem,
which are reported on in this volume.
One can consider searching for 
general methods that apply  to a large variety of related iterations. Such general methods  
as are known give useful
information, and answer some questions about iterates of the $3x+1$ function. 
Nevertheless it is fair to say
that  they do not begin to answer the central question:

\smallskip
{\em What is the 
ultimate fate under iteration of such maps over all time?}

\smallskip
My personal viewpoint is that the $3x+1$ problem is somewhat dangerous, and that  it is prudent
not to focus on resolving the $3x+1$ conjecture as an
immediate goal. Rather, one  might first look for more structure in the problem. 
Also one might profitably 
 view the problem as  a ``test case",  to which one may from time to time apply new results arising from  
the ongoing development of mathematics. When new theories and new methods
are discovered, the $3x+1$  problem may be used as  a testbed to assess their  power,  
 whenever circumstances permit. 


To conclude, let  us remind ourselves, following Hilbert \cite[p. 438]{Hi00}: 

\begin{quotation}
The mathematicians of past centuries were accustomed to devote themselves to the
solution of difficult particular problems with passionate zeal. They 
knew the value of difficult problems.
\end{quotation} 

\smallskip
\noindent The  $3x+1$ problem stands before us as a beautifully simple
question.  It is hard to resist exploring its structure.
We should not exclude it from the mathematical universe just
because we are unhappy with its difficulty.  It is a fascinating and addictive problem.

\bigskip



\paragraph{\bf Acknowledgments.}
I am grateful to  Michael Zieve and Steven J. Miller 
each for detailed readings and corrections.
  Marc Chamberland, Alex Kontorovich, and Keith R. Matthews  also made
  many  helpful comments. The  figures are due to Alex Kontorovich, adapted from
  work on  \cite{KL09}.
I thank Andreas Blass for  useful comments on incompleteness results
 and algebraic structures.
The  author was supported by NSF Grants DMS-0500555 and DMS-0801029.

%
%
%


\begin{thebibliography}{99}
\bibitem{Ga10}
S. Garoufalidis,
The degree of a $q$-holonomic sequence is a quadratic quasi-polynomial,
eprint: {\tt arxiv:1005.4580v1}
\end{thebibliography}
\end{document}